\documentclass{svjour3}                     
\smartqed  
\usepackage{graphicx}
 \usepackage{mathptmx}      
\usepackage{latexsym}
\usepackage[shortlabels]{enumitem}
\usepackage{amsmath, amssymb, amsfonts}
\usepackage{color}
\usepackage[T1]{fontenc}
\usepackage[cp1250]{inputenc}
\usepackage{authblk}
\usepackage{nicematrix}
\usepackage{yfonts}

\newtheorem{thm}{Theorem}
\newtheorem{lem}[thm]{Lemma}
\newtheorem{cor}[thm]{Corollary}

\newtheorem{exm}[thm]{Example}
\newtheorem{rmk}[thm]{Remark}

\begin{document}

\title{On the extension of Batty's theorem on the semigroup asymptotic stability}
\author{Grigory M. Sklyar    \and
       Piotr Polak  \and
        Bartosz Wasilewski 
        }
\institute{G.M. Sklyar (corresponding author) \at
            Institute of Mathematics, University of Szczecin, Wielkopolska 15, 70-451 Szczecin, Poland;
            \\ Faculty of Computer Science and Information Technology, West Pomeranian University of Technology in Szczecin, Zolnierska 49, 71-210 Szczecin, Poland
           \\ \email{Grigorij.sklyar@zut.edu.pl}
           orcid 0000-0003-4588-9926\and
P. Polak \at
            Institute of Mathematics, University of Szczecin, Wielkopolska 15, 70-451 Szczecin, Poland
           \\\email{piotr.polak@usz.edu.pl}
            orcid 0000-0002-8043-3505\and
B. Wasilewski \at
         Institute of Mathematics, University of Szczecin, Wielkopolska 15, 70-451 Szczecin, Poland, Institute of Mathematics, University of Szczecin, Doctoral School, Mickiewicza 16, 70-383 Szczecin, Poland 
         \email{bartosz.wasilewski@phd.usz.edu.pl}
         orcid 0000-0002-9336-5273}


\date{Received: date / Accepted: date}
\maketitle

\begin{abstract}
The well-known Batty's theorem states that if a $C_0$-semigroup $T(t)$ is bounded and the spectrum of the generator $A$ is contained in the open left-half plane of $\mathbb{C}$, then $\|T(t)A^{-1}\|$ tends to $0$. This can be thought of as a particular case of a more general property that, for $\omega_0>-\infty$ and $(\omega_0+i\mathbb{R})\cap \sigma(A)=\emptyset$ it holds $\|T(t)(A-\omega_0 I)^{-1}\|/\|T(t)\|$ tends to 0. We show that it is true for $\|T(t)\|$ regular enough, however we give examples of unbounded semigroups, with the spectrum of the generator not contained in the open left-half plane of $\mathbb{C}$, with the above property. Moreover we give a more general sufficient condition for this property to hold, thus extending Batty's theorem.

\end{abstract}

\section{Introduction}
\label{Introduction}

The asymptotic behavior of semigroups and their orbits has been a subject of an intense  study for the last few decades, see e.g. \cite{BoTo10}, \cite{BaBoTo16}, \cite{Hu93}. 
In \cite{SkSh82}, \cite{LuPh88}, \cite{ArBa88} the authors obtained necessary and sufficient conditions for strong stability of bounded semigroups. In particular, for the case of
$\sigma(A)$, the spectrum of the generator,  is contained in the open left-half plane $\{z\in \mathbb{C}:\Re( z) <0\}$ the semigroup is strongly stable. Due to the Banach-Steinhaus Theorem, if the growth bound $\omega_0(T)=0$ this stability cannot be uniform. However,  due to the works \cite{Ba94}, \cite{BaDu08} of Batty and Batty and Duyckaerts, we have the following theorem
\begin{thm} \label{ThmBa94}
Let $T=\{T(t)\}_{t\geq0}$ be a bounded $C_0$-semigroup acting on a Banach space $X$ and let $A$ be its generator. Then $\|T(t)A^{-1}\| \to 0$ as $t\to +\infty$ if and only if $\sigma(A)\cap (i\mathbb{R})= \emptyset$.
\end{thm}
The above means that for a bounded semigroup $T$ for which
\begin{equation}\label{rez}
    \sigma(A)\subset \{z\in \mathbb{C}:\Re (z) <0\},
\end{equation}
the operator-valued function $T(t):\mathbb{R}^+\ni t\to\mathcal{L}(D(A),X)$
 tends to $0$ as $t\to\infty.$ 
With this being the case, we call the semigroup, after Batty and Duyckaerts, who defined it for the case of bounded semigroups \cite{BaDu08}, as \textit{semi-uniformly stable}.
Moreover, the semi-uniform stability may occur even for unbounded semigroups (see \cite{SkPo17}, for example). 
For the case of unbounded semigroups  it was shown in \cite{Sk15} that the condition \eqref{rez} remains necessary for $\|T(t)A^{-1}\| \to 0$. On the other hand for an unbounded semigroup $T$ with $\omega_0(T)\geq0$ the concept of semi-uniform stability led us to consider a more general property: 
\begin{equation}\label{Teza}
\lim_{t\to+\infty}\frac{\|T(t)R_\mu\|}{\|T(t)\|}=  0, \quad  \textnormal{ for }  \mu \not\in \sigma(A).
\end{equation}
where by $R_\mu$ we mean the resolvent of the semigroup generator at the point $\mu\notin \sigma(A)$. This property can be thought of as the growth rate of the semigroup truncated to the domain of the generator being slower than the growth rate of the semigroup.
For bounded semigroups with $\omega_0(T)=0$, \eqref{Teza} clearly reduces  to the semi-uniform stability. A question arises here about the necessary and sufficient conditions for the property \eqref{Teza} to hold for general semigroups. The condition
\begin{equation}\label{warunek}
    (\omega_0(T)+i\mathbb{R})\cap\sigma(A) = \emptyset
\end{equation} is not necessary for \eqref{Teza} to occur, as it is shown in Example \ref{ex1}. In this example the behavior of $\|T(t)\|$ and $\|T(t)R_\mu\|$ is easy to predict due to existence of an orthonormal basis. Similarly, \eqref{Teza} can be verified for a $C_0$-group $\{T(t)\}_{t\in\mathbb{R}}$ such that the spectrum of the generator $A$ is discrete and the eigenvalues are uniformly separated. Indeed, in this case, due to \cite{XuYu05}, \cite{Zw10}, the corresponding eigenvectors constitute a Riesz basis and the problem of verifying \eqref{Teza} may be reduced to solving the problem in the invariant subspaces. In general Banach spaces this problem becomes more complicated. The main goal of this paper was to present a  sufficient condition for \eqref{Teza} to hold, in the case when eigenvectors do not necessarily constitute a Riesz basis. In section \ref{application} we show that the property \eqref{Teza} holds for a class of unbounded semigroups for which $\sigma(A)\subset (i\mathbb{R})$, $\sigma(A)$ is countable, and consists of simple eigenvalues only.
\section{Main result}\label{main}

First we show that the condition \eqref{warunek} is not necessary for the property \eqref{Teza} to hold. Below we give an example of an unbounded $C_0$-semigroup with $\omega_0=0$ for which \eqref{Teza} holds, despite the fact that $\sigma(A)\cap (i\mathbb{R}) \neq \emptyset$.
\begin{exm}\label{ex1}
Consider a separable Hilbert space $H$ with the orthonormal basis $\{e_n\}_{n\in\mathbb{N}}$ and put 
\begin{equation*}
    T(t)e_0=e^{it}e_0,\quad T(t)e_{2k-1}=e^{(ik-\frac 1k)t}e_{2k-1}, \quad T(t)e_{2k}=e^{(ik-\frac 1k)t}(te_{2k-1}+e_{2k}), 
\end{equation*}
for $k=1,2,\ldots$ The above defines a $C_0$-semigroup $T=\{T(t)\}_{t\geq0}$ on $H$.  It is easy to see that on the invariant subspace \begin{equation*}
    H_1={\rm span}\{e_0\},
\end{equation*}the operators $T(t)$ and $T(t)R_\mu$ are uniformly bounded for $t\geq 0$.
It is less obvious that on the complementary subspace
\begin{equation*}
    H_2=\overline{{\rm span}\{e_1,e_2,\ldots\}}, 
\end{equation*} the norm of the semigroup behaves as following:
\begin{equation*}
    \|T(t)\|\sim t.
\end{equation*}
 In particular, this implies $\omega_0 = 0.$ Also, direct computations (or applying the result from \textnormal{\cite{SkPo17}}) show that
\begin{equation*}
    \|T(t)R_\mu\|\leq M, \quad t\geq 0.
\end{equation*} This means that \eqref{Teza} holds despite
\begin{equation*}
    \{i\}\subset\sigma(A)\cap(i\mathbb{R})\neq \emptyset.
\end{equation*}
\end{exm}
However, we will prove that the condition \eqref{warunek} is  sufficient for \eqref{Teza} to hold. This will follow from the next theorem which is the main result of this work:
\begin{thm}\label{main} Let $T = \{T(t)\}_{t\geq 0}$ be a semigroup on a Banach space $X$, not necessarily bounded, with the growth bound $\omega_0>-\infty$ and the generator $A$. Suppose $f(t):\mathbb{R}^+\to\mathbb{R}^+$ is a positive function with concave downwards logarithm $\log(f(t))$ which approximates the semigroup norm $\|T(t)\|$ in the following sense
    \begin{equation}{\label{a}}
        \|T(t)\|\leq f(t),\quad t\geq 0,
    \end{equation}
     \begin{equation}\label{b}
 \limsup_{t\to+\infty} \,\frac{ \|T(t)\|}{f(t)}  = a>0.
   \end{equation}
   
Assume further that 
\begin{enumerate}[(a)]
    \item for any $\lambda \in \sigma(A) \cap (\omega_0+i\mathbb{R})$ there exists a regular bounded curve $\Gamma_\lambda$ enclosing $\lambda$, such that $\Gamma_\lambda \cap \sigma(A) = \emptyset;$ 
    \item for any $\lambda \in \sigma(A) \cap (\omega_0+i\mathbb{R})$
    \begin{equation}\label{Riesz}
        \underset{t \to +\infty}{\lim} \frac{\| T(t)P_{\Gamma_\lambda} \|}{f(t)}= 0,
    \end{equation}
    \end{enumerate}
    where $ P_{{\Gamma}_\lambda}$ is the Riesz projection associated with the curve $P_{{\Gamma}_\lambda}$. 
Then
\begin{equation}{\label{main1}}
     \underset{t \to +\infty}{\lim} \frac{\|T(t)R_\mu\|}{f(t)}=0,
    \end{equation}
for fixed $\mu \not\in \sigma(A)$.
\end{thm}
Before the proof of the theorem a few remarks are in order:
\begin{itemize}
\item a constructive proof of existence of such a function $f$ satisfying \eqref{a} and \eqref{b} for an arbitrary semigroup is given in \textnormal{\cite{Sk15}}; 
\item without loss of generality we only prove the Theorem \eqref{main} in the case of $\omega_0=0$. Indeed, for arbitrary $\omega_0$ one can consider the shifted semigroup $\{e^{-\omega_0t}T(t)\}_{t\geq 0}$; 
\item we will clarify the connection between \eqref{main1} and \eqref{Teza} at the end of the proof.
\end{itemize}
In the proof we will use the construction of the special operator-valued semigroup introduced in \cite{Sk15}.  
Let $\widetilde{X}\subset \mathcal{L}(X)$ be defined as 
\begin{equation}
    \widetilde{X} = \overline{\{DR_{\mu}(A),\quad D\in \mathcal{L}(X)\}},\quad \mu \not\in\sigma(A),
\end{equation}
 where $\overline{Q} $ denotes the closure of the linear set $Q$ (with respect to the operator norm).  Since $\widetilde{X}$ is a closed subspace of a Banach space  $\mathcal{L}(X)$, it also is a Banach space. It is clear that $\widetilde{X}$ does not depend on the choice of $\mu.$ For the given semigroup \{$T(t)\}_{t\geq 0 }$ on  the space $X$, let us introduce a semigroup on the space $\widetilde{X}$ by:
\begin{equation}{\label{TB}}
    \widetilde{T}(t)\widetilde{B} = \widetilde{B}T(t), \quad \widetilde{B}\in \widetilde{X}, \quad t\geq 0.
\end{equation}
 Important properties of this semigroup were shown in \cite{Sk15}, namely that $\{\widetilde{T}(t)\}_{t\geq 0}$ forms a $C_0$-semigroup on $\widetilde{X}$, and that 
\begin{itemize}
  \item for $A$ and $\widetilde{A}$ being the generators of $\{{T}(t)\}_{t\geq 0}$ and $\{\widetilde{T}(t)\}_{t\geq 0}$, respectively, it holds that \begin{equation}\label{aSK}
      \sigma(\widetilde{A})\subset \sigma(A);
  \end{equation}
  \item for $\widetilde{B}\in \widetilde{X}$ and $\mu \not\in \sigma(A)$, it holds that \begin{equation}\label{bSk}
      (\widetilde{A} - \mu  I)^{-1}\widetilde{B} = \widetilde{B}(A-\mu I)^{-1}.
  \end{equation} 
\end{itemize}
We will also use the following lemma 
\begin{lem}{\label{LemEN}}\textnormal{\cite{LuPh88}}
Let $\{T(t)\}_{t\geq 0}$ be a strongly continuous semigroup of isometries on a Banach space $X$ and denote its generator by $A.$ Then one of the following two cases holds
\begin{itemize}
    \item $\sigma(A) = \{\mu\in\mathbb{C}:\Re (\mu)\leq 0\}$;
    \item $\sigma(A) \subset(i\mathbb{R})$ and the above semigroup extends to a strongly continuous group of isometries.
    \end{itemize}
\end{lem}
Note that Lemma \ref{LemEN} implies that, for a semigroup of isometries, if $\partial(\sigma(A))\neq (i\mathbb{R})$, then $\sigma(A) =  \partial(\sigma(A))\varsubsetneq (i\mathbb{R}), $ where $\partial $ denotes the boundary of a set.
The proof of Theorem \ref{main} is based on the idea used in \cite{SkPo19}.
\\ \textit{Proof of Theorem \ref{main}}.
 
 \noindent \\ Assume that (\ref{main1}) does not hold, which means that 
\begin{equation}\label{As}
     0\neq \underset{t \to +\infty}{\lim\sup} \frac{\|T(t)R_\mu\|}{f(t)}=
      \underset{t \to +\infty}{\lim\sup} \frac{\|R_\mu T(t)\|}{f(t)}=
     \underset{t \to +\infty}{\lim\sup} \frac{\|\widetilde{T}(t)R_\mu\|}{f(t)}.
\end{equation}
Let us define a following seminorm on $\widetilde{X}$:
\begin{equation*}
    l(\widetilde{B}) = \underset{t \to +\infty}{\lim \sup}\frac{\|\widetilde{T}(t)\widetilde{B}\|}{f(t)},\quad \widetilde{B}\in \widetilde{X}.
\end{equation*}
It follows from \eqref{As} that the quotient space $\widetilde{X}\slash \ker{l} = \{\widehat{B} = \widetilde{B}+\ker l:\widetilde{B}\in \widetilde{X}\}$ is non-zero. This space can be equipped with a norm different from the natural one  ($\|\widehat{B}\|_N:=\inf\{\|\widetilde{B}\|:\widetilde{B}\in \widehat{B}\}$) of the following form 
\begin{equation*}
    \|\widehat{B}\|':=l(\widetilde{B}), \quad \widetilde{B}\in\widetilde{X}.
\end{equation*}
Note that, since $\|\widetilde{T}(t)\|\leq\|T(t)\|\leq f(t)$ (see (\ref{TB}), (\ref{a})), for all $\widetilde{B}\in \widetilde{X}$, 
\begin{equation*}
     l(\widetilde{B}) = \underset{t \to +\infty}{\lim \sup}\frac{\|\widetilde{T}(t)\widetilde{B}\|}{f(t)}\leq \|\widetilde{B}\|
\end{equation*}
holds, which means that $\|\widehat{B}\|'\leq \|\widehat{B}\|_N$ and the space $(\widetilde{X}\slash \ker{l}, \|\cdot\|')$ may be incomplete. Its completion w.r.t. the norm $\|\cdot\|' $ is denoted by $\widehat{X}$. Let us define the family of operators $\widehat{T}(t), t\geq 0$ by the formula \begin{equation*}
    \widehat{T}(t)\widehat{B} = \widetilde{T}(t)\widetilde{B} +\ker l,\quad \widehat{B}\in \widetilde{X}\slash \ker l\subset \widehat{X}.
\end{equation*}
We will now prove that $\widehat{T}(t), t\geq 0 $ is a family of isometries on $\widetilde{X}\slash \ker{l}$, w.r.t. the norm $\|\cdot\|'$. This follows from an assertion from real analysis. Namely, \cite{Taiw}, \cite{Dok}
let $h(t)$ be a real non-negative function defined on the positive semi-axis $\mathbb{R}^+ = \{t:t\geq0\}$ and such that 
\begin{itemize}
    \item $\forall \epsilon>0\textnormal{ } \exists C_\epsilon,  \textnormal{ s.t. }$  $h(t) \leq C_\epsilon +\epsilon t,\quad t\geq0;$
    \item$h(t)$ is concave downwards.
\end{itemize}
Then for any $t>0$ the following holds:
\begin{equation*}
    \underset{s \to +\infty}{\lim}(h(t+s)-h(s)) = 0.
\end{equation*}
\begin{flushleft}
By applying this assertion to $f$ meeting conditions of Theorem \ref{main}, we obtain
\end{flushleft}
\begin{equation}
  \underset{s \to +\infty}{\lim}e^{\log(f(t+s))-\log(f(s))} =   \underset{s \to +\infty}{\lim}\frac{f(t+s)}{f(s)} = 1.
\end{equation}
Now, applying this result, we get
\begin{equation*}
    \|\widehat{T}(t)\widehat{B}\|' = \underset{s \to +\infty}{\lim \sup}\frac{\|\widetilde{T}(t+s)\widetilde{B}\|}{f(t+s)}\frac{f(t+s)}{f(s)} = \|\widehat{B}\|',\quad \textnormal{ for } \widehat{B}\in \widetilde{X}\slash \ker{l}.
\end{equation*}
Thus, $\widehat{T}(t), t\geq 0 $ is a family of isometries on $\widetilde{X}\slash \ker{l}$, w.r.t. the norm $\|\cdot\|'$. $ $ It is easy to check that for each $t\geq 0,$ $\widehat{T}(t)$ extends to an isometry on $\widehat{X}$ and the family $\widehat{T}(t), t\geq 0 $ is a $C_0$-semigroup of isometries. Moreover, one can check that  
\begin{equation}\label{RGen}
\widehat{A}\widehat{B} = \widetilde{A} \widetilde{B}+\ker l,  \textnormal{ and} 
    \end{equation}
    \begin{equation*}
         R(\widehat{A}, \mu)\widehat{B} = R(\widetilde{A},\mu)\widetilde{B}+\ker l
    \end{equation*}
for $\widehat{B}\in \widetilde{X}$, where $\widetilde{A}$ and $\widehat{A}$ are generators of $\{\widetilde{T}(t)\}_{t\geq 0}$ and $\{\widehat{T}(t)\}_{t\geq 0}$, respectively and $ R(\widetilde{A}, \mu)$ and $ R(\widehat{A}, \mu)$ are the respective resolvent operators at the point $\mu$. It follows from assumption (a) of Theorem \ref{main} and \eqref{aSK} that
\begin{align}\label{spec_1}
   & (i\mathbb{R})\not\subset\sigma(A)  \\ 
   & (i\mathbb{R})\not\subset\sigma(\widetilde{A}). \nonumber
\end{align}
On the other hand, it is shown in \cite{Taiw},\cite{Dok} that  \begin{equation*}
    \partial(\sigma(\widehat{A}))\cap(i\mathbb{R}) \subset  \sigma(\widetilde{A})\cap(i\mathbb{R}),
\end{equation*}where $\partial$ denotes the boundary of a set. This, along with Lemma \ref{LemEN} and \eqref{spec_1}, implies that 
\begin{equation}{\label{spec}}
\partial\sigma(\widehat{A})=\sigma(\widehat{A})\subset\sigma(\widetilde{A})\cap(i\mathbb{R})\neq (i\mathbb{R}).
\end{equation} 
Therefore, again due to Lemma \ref{LemEN}, $\{\widehat{T}(t)\}_{t\geq 0}$ extends to a $C_0$-group of isometries. Now, since $\widehat{A}$ is a generator of a $C_0$-group of isometries, its spectrum has to be non-empty (see e.g. \cite{Neerven})
\begin{equation*}
    \sigma(\widehat{A})\neq\emptyset.
\end{equation*}
 By combining the above with (\ref{spec}) and (\ref{aSK}), we obtain:
\begin{equation}\label{specinc2}
 \emptyset\neq   \sigma(\widehat{A})\subset \sigma(\widetilde{A})\cap(i\mathbb{R})\subset\sigma(A)\cap(i\mathbb{R).}
\end{equation}
Note that in the case $\sigma(A)\cap(i\mathbb{R})=\emptyset$ we obtain here a contradiction. This means that
\begin{equation*}
    \underset{t \to +\infty}{\lim} \frac{\|T(t)R_\mu\|}{f(t)}=0.
\end{equation*}
Now, for the case when $\sigma(A)\cap(i\mathbb{R})\neq\emptyset$,  let us fix $\lambda$ such that
\begin{equation*}
    \lambda\in \sigma(\widehat{A})\subset \sigma(A)\cap(i\mathbb{R}).
\end{equation*} 
It follows from the assumption (a), \eqref{specinc2}, and \eqref{aSK} that there exists a bounded curve $\Gamma_{\lambda}$ enclosing $\lambda$, such that
\begin{equation*}
    \Gamma_{\lambda}\cap\sigma(\widehat{A}) =  \Gamma_{\lambda}\cap\sigma(\widetilde{A})= \Gamma_{\lambda}\cap\sigma(A) = \emptyset.
\end{equation*} Let $\widetilde{P}_{\Gamma_{\lambda}}$ and $\widehat{P}_{\Gamma_{\lambda}}$ be the
Riesz projections in $\widetilde{X}$ and $\widehat{X}$, respectively, corresponding to the curve 
$\Gamma_{\lambda}$. One can see from \eqref{RGen}, that for $\widehat{B} \in \widetilde{X}/\ker l$
\begin{equation}\label{Riesz2}
    \widehat{P}_{\Gamma_{\lambda}}\widehat{B} = \widetilde{P}_{\Gamma_{\lambda}} \widetilde{B}
+ \ker l.
\end{equation}
   Furthermore, the projections $\widetilde{P}_{\Gamma_{\lambda}}$ and $\widehat{P}_{\Gamma_{\lambda}}$ split the spaces $\widetilde{X}$ and $\widehat{X}$ into direct sums $\widetilde{Z}_1+\widetilde{Z}_2$ and $\widehat{Z}_1$ + $\widehat{Z}_2$, respectively, so that \begin{align*}
    & \widetilde{Z}_1 :=  \widetilde{P}_{\Gamma_{\lambda}} \widetilde{X},\\
    & \widetilde{Z}_2 :=  (I-\widetilde{P}_{\Gamma_{\lambda}}) \widetilde{X},\\
     & \widehat{Z}_1 :=  \widehat{P}_{\Gamma_{\lambda}} \widehat{X},\\
    & \widehat{Z}_2 :=  (I-\widehat{P}_{\Gamma_{\lambda}}) \widehat{X}.
  \end{align*}
  Clearly the spectra of the restricted operators $\widetilde{A}|_{\widetilde{Z}_1}$ and $\widetilde{A}|_{\widetilde{Z}_2}$ are intersections of $\sigma(\widetilde{A})$ with regions inside and outside $\Gamma_{\lambda}$, respectively, with an analogous property for $\sigma(\widehat{A})$. Now, since the set $\sigma(\widehat{A})$ is a boundary set, it consists only of approximate eigenvalues (see e.g. [\cite{EnNa00} VI, prop. 1.10]). This means that
for the chosen $\lambda$ there exists a sequence $\{\widehat{B}_k\}:\|\widehat{B}_k\|'=1$ such that
\begin{equation}{\label{5}}
    \|\widehat{A}\widehat{B}_k-\lambda\widehat{B}_k\|'\to 0 \textnormal{ as } k\to \infty.
\end{equation}
Now, $\{\widehat{B}_k\}$ can be split into a sequence
\begin{equation*}
    \widehat{B}_k = \widehat{B}_{k}^{(1)}+\widehat{B}_{k}^{(2)},
\end{equation*}where
    \begin{equation*}
        \widehat{B}_{k}^{(1)}\in \widehat{Z}_1,\quad \widehat{B}_{k}^{(2)}\in \widehat{Z}_2.
    \end{equation*}
Then it follows from (\ref{5}), that
\begin{align*}
    & \|\widehat{A}\widehat{B}_k^{(1)}-\lambda\widehat{B}_k^{(1)}\|'\to 0,\\
    & \|\widehat{A}\widehat{B}_k^{(2)}-\lambda\widehat{B}_k^{(2)}\|'\to 0,
\end{align*}
as $k\to\infty$. Subsequently,
\begin{equation*}
\|    \widehat{B}_k^{(2)}\|\to 0,
\end{equation*}
since otherwise $\lambda$ would belong to $\sigma(\widehat{A}|_{\widehat{Z}_2})$, giving a contradiction. In consequence 
\begin{equation*}
\|\widehat{B}_{k}^{(1)}\|'\geq\frac{1}{2}    
\end{equation*}  for $k$ large enough. Furthermore, by the density of ${\widetilde{X}\slash \ker l}$ in $\widehat{X}$ and by the boundedness of $\widehat{A}|_{Z_1}$, $\widehat{B}_{k}^{(1)}$ can be chosen from $\widehat{P}_{\Gamma_{\lambda}}(\widetilde{X}\slash \ker l)\subset\widehat{Z}_1$. Subsequently, from \eqref{Riesz2}, we get
\begin{equation*}
    \widehat{B}_{k}^{(1)} = \widehat{P}_{\Gamma_{\lambda}}\widehat{B}_{k}=\widetilde{P}_{\Gamma_{\lambda}}{\widetilde{B}}_{k} + \ker l,
\end{equation*}
for some sequence $\widetilde{B}_{k}\in \widetilde{X}$. Then the following estimate holds
\begin{equation}{\label{last}}
    \begin{split}
    \frac{1}{2}	\leq \|\widehat{B}_{k}^{(1)}\|' & =\|\widehat{P}_{\Gamma_{\lambda}}\widehat{B}_{k}\|'= \|\widetilde{P}_{\Gamma_{\lambda}}{\widetilde{B}}_{k} + \ker l\|' = l(\widetilde{P}_{\Gamma_{\lambda}}{\widetilde{B}_{k}})\\ & = \underset{t \to +\infty}{\lim \sup}\frac{\|\widetilde{T}(t)\widetilde{P}_{\Gamma_{\lambda}}\widetilde{B}_{k}\|}{f(t)},
    \end{split}
\end{equation}
for $k$ large enough. By integrating the equation (\ref{bSk}), we obtain
\begin{equation}{\label{p}}
    \widetilde{P}_{\Gamma_{\lambda}}\widetilde{B}_k= \int_{\Gamma_{\lambda}}  (\widetilde{A} - \mu  I)^{-1}\widetilde{B}_k d\mu = \int_{\Gamma_{\lambda}}\widetilde{B}_k(A-\mu I)^{-1}d\mu=  \widetilde{B}_k P_{\Gamma_{\lambda}}.
\end{equation} 
Where we have used the analicity of the 
resolvent operator function and the boundedness of $\widetilde{B}_k$ as an operator from $\mathcal{L}(X)$ to $\mathcal{L}(X)$ (treated as a multiplication operator).  Now recall that we assumed
\begin{equation*}\tag{6}
        \underset{t \to +\infty}{\lim} \frac{\| T(t)P_{\Gamma_\lambda} \|}{f(t)}= 0.
    \end{equation*}
 Using (\ref{last}), (\ref{p}), \eqref{Riesz}, and the definition of $\{\widetilde{T}(t)\}_{t\geq 0 }$ (see \eqref{TB}), we get

\begin{gather*}
  \frac{1}{2}\leq \|\widehat{B}_{k}^{(1)}\|'= l(\widetilde{P}_{\Gamma_{\lambda}}{\widetilde{
  B}_{k}}) = \underset{t \to +\infty}{\lim \sup}\frac{\|\widetilde{T}(t)\widetilde{P}_{\Gamma_{\lambda}}{\widetilde{B}}_{k}\|}{f(t)}= \underset{t \to +\infty}{\lim \sup}\frac{\|{\widetilde{B}}_{k}P_{\Gamma_{\lambda}} T(t)\|}{f(t)} \\\leq
 \underset{t \to +\infty}{\lim \sup}\frac{\|{\widetilde{B}}_{k}\|\|P_{\Gamma_{\lambda}} T(t)\|}{f(t)} =0. 
\end{gather*}
This yields a contradiction, thus
 \begin{equation*}
     \underset{t \to +\infty}{\lim} \frac{\|\widetilde{T}(t)R_\mu\|}{f(t)}=\underset{t \to +\infty}{\lim} \frac{\|T(t)R_\mu\|}{f(t)}=0.
 \end{equation*}
\begin{flushright}
$\square$
\end{flushright}
\begin{rmk}\label{remark}
 For bounded semigroups ($\|T(t)\|\leq M$ for $t\geq0$) with  $\sigma(A)\cap (i\mathbb{R}) = \emptyset$ and $\omega_0(T) = 0$, by taking $f\equiv M$ one can easily see that Theorem \ref{main} implies the sufficiency part in Theorem \ref{ThmBa94}.
Also note that, if
$\|T(t)\|$ has a concave downwards logarithm or
\begin{equation*}
    \|T(t)\|\sim f(t),
\end{equation*}
for some f(t) with concave downwards logarithm (e.g. $f(t) = t^\alpha e^{\beta t}$, for some $\alpha \in \mathbb{R}^+, \beta \in \mathbb{R}$), 
 then the assertion of
Theorem \ref{main} takes the following form
\begin{equation*}
\frac{\|T(t)R_\mu\|}{\|T(t)\|}\to 0, \quad \textnormal{ as } t\to\infty.
\end{equation*}
Where by $\|T(t)\|\sim f(t)$ we mean that there exists $c,C>0$, such that 
\begin{equation*}
 cf(t)\leq \|T(t)\|\leq Cf(t), \quad t\geq t_0   
\end{equation*} for some $t_0\geq0$.
\end{rmk}

\section{Generator with countable pure imaginary simple spectrum}\label{application} Here we give some examples of semigroups with generators having a countable pure imaginary simple spectrum for which our result can be applied. Let us consider the case when the eigenvalues are uniformly separated, i.e., $\inf \{ |\lambda_k-\lambda_m|:k,m\in\mathbb{N}, k\neq m\}>0 $ and the eigenvectors are linearly dense. Then, due to the works of \cite{XuYu05} and \cite{Zw10} we know that the eigenvectors form a Riesz basis in $H$. It follows the semigroup is bounded and due to Batty's and Duyckaerts' theorem (\cite{BaDu08}), $\|T(t)R_\mu\|$ $\not\to$ 0 as $t\to\infty$, hence \eqref{Teza} cannot hold.  
Consider the following examples for which the eigenvalues are no longer uniformly separated, thus  allowing the semigroup to be possibly unbounded. 

\begin{exm}
Let $\{e_n\}_{n=1}^{\infty}$ be the orthonormal basis of a Hilbert space $H$. Define the operator $A:D(A)\subset H\to H$ as follows:
\begin{equation*}
A|_{H_n} := A_n := 
\begin{bmatrix}
ni+\frac{i}{n} & 1  \\
0 &  ni-\frac{i}{n}
\end{bmatrix},
\end{equation*}
where $H_n= {\rm span} \{e_{2n-2},e_{2n-1}\},\textnormal{ } n=2,3,4\ldots$. 
For each $n\geq 2$ consider the curve $\Gamma_n$ enclosing the pair of eigenvalues $i(n+\frac1n), i(n-\frac1n)$, then the image of the Riesz projection corresponding to the curve is $H_n$. 
One can directly check that

\begin{equation*}
e^{A_nt} := T_n(t) = e^{tni} 
\begin{bmatrix}
e^{i\frac tn} &\textnormal{ } n\sin\frac tn \\[2mm]
0 &  e^{-i\frac tn}
\end{bmatrix},
\end{equation*}
Since $\|T(t)\|=\displaystyle\sup_{n\geq2}\left\|T_n(t)\right\|$, we have
\begin{equation*}\label{eAt}
\|T(t)\| \sim t.
\end{equation*}
It is easy to see, that $f(t) :=t$ has a concave downwards logarithm and has the desired properties \eqref{a}, \eqref{b}. Clearly assumptions (a) and (b) of Theorem \ref{main} are satisfied. Therefore \eqref{Teza} holds, i.e, 
\begin{equation}\label{limit1}
    \frac{\|T(t)A^{-1}\|}{t}\to 0,\quad t \to \infty.
\end{equation}
Moreover, for this simple case, we can calculate the
decay rate of \eqref{limit1}, namely
\begin{equation*}
T_n(t)A_n^{-1} = \frac{in}{1-n^4}e^{tni}
\begin{bmatrix}
(n^2-1) e^{i\frac{t}{n}} & \textnormal{ }(n^2-1)n\sin\frac tn +ine^{-i\frac{t}{n}}  \\ 
0 &  (n^2 + 1)e^{-i\frac{t}{n}}
\end{bmatrix},
\end{equation*}
hence
\begin{equation*}
\|T(t)A^{-1}\|=\sup_{n\geq2}\|T_n(t)A_{n}^{-1}\|\sim 1, \quad t\geq0.
\end{equation*}
 Finally, it follows that
\begin{equation*}
\frac{\left\|T(t)A^{-1}\right\|}{\|T(t)\|}\sim \frac{1}{t} \to 0,\quad t \to \infty.
\end{equation*}
\end{exm}

 Below we give an example of a family of unbounded semigroups that have simple countable purely imaginary spectrum and the eigenvectors are linearly dense but do not form a Riesz basis. This family was  described in \cite{SMP}  and \cite{SM}. The elements of this family are constructed as follows. Let $(H, \|\cdot\|)$ be a Hilbert space with the orthonormal basis $\{e_n\}_{n=2}^{\infty}$. For the sequence
\begin{equation*}
    \lambda_n = i\log n,\, n=2,3,\ldots.
\end{equation*}
 define a semigroup $T $ by 
 \begin{equation*}
     T(t)e_n = e^{t\lambda_n}e_n,
 \end{equation*}
 For a given $N\in\mathbb{N}/\{0\}$ we are able to choose a new norm $\|\cdot\|_N$ on $H$, dominated by $\|\cdot\|$ such that:
\begin{itemize}
    \item $T$ naturally extends to a $C_0$-semigroup $\widetilde{T}$ on the completion of $(H,\|\cdot\|_N)$, say $\widetilde{H}_N$;
    \item there exist constants $m,M>0$ such that 
\begin{equation}\label{Hnorm}
    mt^N\leq \|\widetilde{T}(t)\|\leq Mt^N + 1, \quad t\geq0.
\end{equation}
\end{itemize}
See \cite{SMP} and \cite{SM} for a detailed construction and estimations. Denote the generator of $\widetilde{T}$ by $\widetilde{A}$. Then \cite{SM}
\begin{equation*}
    \sigma(\widetilde{A}) =\sigma_P(\widetilde{A})= \underset{n\geq 2}{\bigcup}i\log n.
\end{equation*}
We are going to show that the semigroup $\widetilde{T}$ meets the assumptions of Theorem~\ref{main}, however first we will compute the $\|\widetilde{T}(t)\widetilde{A}^{-1}\|/{\|\widetilde{T}(t)\|}$ "by hand" for the case of $N = 1$. Before we do that we should show some basic properties of the space  $(H,\|\cdot\|_N)$, as shown in \cite{SMP}, \cite{SM}. Consider the backward difference operator 
\begin{equation*}
    \Delta = 
\begin{bNiceMatrix}
1 & 0 & 0 & 0 & \textnormal{ }\cdots\\
-1 & 1 & 0 & 0 & \cdots\\
0 & -1 & 1 & 0 & \cdots\\
0 & 0 & -1 & 1 & \cdots\\
\vdots & \vdots & \vdots & \vdots & \ddots
\end{bNiceMatrix}.
\end{equation*}
The space $(H,\|\cdot\|_N)$ is defined as the completion of
\begin{equation*}
    \Big \{ x  =(\textgoth{f})\sum^{\infty}_{n = 2} c_n e_n : \{c_n\}_{n=1}^\infty \in l_2(\Delta^N)\Big\},
\end{equation*}
with respect to the norm on this space defined as:
\begin{equation}\label{norm}
    \|x\|_N = \Big \|(\textgoth{f})\sum^{\infty}_{n = 2} c_n e_n\Big\|_N = \Big\|\sum^{\infty}_{n = 2} \sum^{N}_{j = 0}(-1)^j C^j_Nc_{n-j}e_n\Big\|,
\end{equation}
where $l_2(\Delta^N) = \big \{ x =\{c_n\}_{n=2}^\infty , c_n \in \mathbb{C}: \Delta^N x\in l_2\big\}$ and $(\textgoth{f})$ denotes the formal series, where the last norm without the subscript denotes the norm in the initial Hilbert space $H$, and $C^j_N$ denote the binomial coefficients $\binom{j}{N}$. The action of the generator, resolvent at the point $0$ and product of the semigroup and the resolvent are as follows:

\begin{equation*}
    \widetilde{A} e_n = i\log(n)e_n, \quad n\geq 2, 
\end{equation*}
\begin{equation*}
    \widetilde{A}^{-1} e_n = \frac{1}{i\log(n)}e_n  \quad n\geq 2,
\end{equation*}
\begin{equation*}
   \widetilde{T}(t) \widetilde{A}^{-1} e_n = \frac{e^{it\log(n)}}{i\log(n)}e_n, \quad n\geq 2,
\end{equation*}
\begin{equation*}
   \widetilde{T}(t) \widetilde{A}^{-1} x = \sum^{\infty}_{n = 2} c_n\frac{e^{it\log(n)}}{i\log(n)}e_n.
\end{equation*}
Let us consider the simplest case of $\widetilde{T}$ when $N=1$.
\begin{exm}
Consider $\widetilde{T}: \widetilde{H}_1 \to \widetilde{H}_1,$ then 
\begin{equation*}
    \|\widetilde{x}\|_1 = \Big(\sum^{\infty}_{n = 2} \Big|c_{n+1} - c_n\Big|^2 + |c_2|^2\Big)^{\frac{1}{2}}, \quad \widetilde{x} \in \widetilde{H}_1.
\end{equation*}
We will prove that, for this case 
\begin{equation}\label{equiv3}
    \frac{ \| \widetilde{T}(t) \widetilde{R}_\mu\|}{\|\widetilde{T}(t)\|}\
   \sim \frac{1}{\log(t)}.
\end{equation}
In further considerations, we will use the following inequality
\begin{equation}\label{Har2}
\sum_{n=1}^{\infty}\frac{|c_{n}|^2}{n^2}\leq 4\sum_{n=1}^{\infty}|c_{n+1}-c_{n}|^2,\quad \{c_n\}_{n=1}^\infty\subset\mathbb{C},
\end{equation}
which is a special case of Hardy's inequality: 
\begin{equation*}
\sum_{n=1}^{\infty}\Big(\frac{1}{n} \sum^{n}_{k=1}a_k\Big)^p\leq \Big(\frac{p}{p-1}\Big)^p\sum_{n=1}^{\infty}a_n^p,\quad a_n\geq 0,
\end{equation*}
for p = 2. To prove \eqref{equiv3} we will estimate $\| \widetilde{T}(t) \widetilde{A}^{-1} \widetilde{x}\|^2_1$. It is given by
\begin{gather*}
  \| \widetilde{T}(t) \widetilde{A}^{-1} \widetilde{x}\|_1^2 =  \sum^{\infty}_{n = 2} \Big|c_{n+1}\frac{e^{it\log(n+1)}}{i\log(n+1)}- c_n\frac{e^{it\log(n)}}{i\log(n)}\Big|^2 + |c_2|^2 \\\leq
    2  \sum^{\infty}_{n = 2} \Big|c_{n+1}\frac{e^{it\log(n+1)}}{i\log(n+1)} -c_{n+1}\frac{e^{it\log(n)}}{i\log(n)}\Big|^2+ 2\sum^{\infty}_{n = 2} \Big| (c_{n+1}-c_n)\frac{e^{it\log(n)}}{i\log(n)}\Big|^2  +|c_2|^2.
\end{gather*}
The second and third elements of the right-hand side of the inequality are clearly bounded by $B\Big(\frac{t}{\log(t)}\Big)^2\|\widetilde{x}\|_1^2$ and $C\Big(\frac{t}{\log(t)}\Big)^2\|\widetilde{x}\|_1^2,$ $B,C>0$, for $t>e$. We only need to look at the first sum then. 
\begin{gather*}
  \sum^{\infty}_{n = 2} \Big|c_{n+1}\frac{e^{it\log(n+1)}}{i\log(n+1)} -c_{n+1}\frac{e^{it\log(n)}}{i\log(n)}\Big|^2 \\= \sum^{\infty}_{n = 2} \Big|\frac{c_{n+1}}{n}\frac{n(e^{it\log(n)}\log(n+1) - e^{it\log(n)}\log(n))}{\log(n+1)\log(n)}  \\+\frac{c_{n+1}}{n}\frac{n(e^{it\log(n)}\log(n) - e^{it\log(n+1)}\log(n))}{\log(n+1)\log(n)} \Big|^2  \\\leq
  2\sum^{\infty}_{n = 2}  \Big|\frac{c_{n+1}}{n}\frac{n(e^{it\log(n)}\log(n+1) - e^{it\log(n)}\log(n))}{\log(n+1)\log(n)}\Big|^2 \\+ \Big| \frac{c_{n+1}}{n}\frac{n(e^{it\log(n)}\log(n) - e^{it\log(n+1)}\log(n))}{\log(n+1)\log(n)} \Big|^2  \\= 2\sum^{\infty}_{n = 2}  \Big|\frac{c_{n+1}}{n}\frac{n\log(1+\frac{1}{n})}{\log(n+1)\log(n)}\Big|^2 +2\sum^{\infty}_{n = 2} \Big|\frac{c_{n+1}}{n}\frac{n(1 - e^{it\log(1+\frac{1}{n})})}{\log(n+1)}\Big|^2.
\end{gather*}
The first of the above sums, due to Hardy's inequality (see \eqref{Har2}), is bounded by $ D\|\widetilde{x}\|^2$, and thus by $ D\Big(\frac{t}{\log(t)}\Big)^2\|\widetilde{x}\|^2$ for $t>e$. We estimate the remaining sum by splitting it into two $t$-dependent sums.
\begin{gather*}
    \sum^{\infty}_{n = 2}  \Big|\frac{c_{n+1}}{n}\frac{n(1 - e^{it\log(1+\frac{1}{n})})}{\log(n+1)}\Big|^2 \\= \sum^{ }_{2\leq n < t}  \Big|\frac{c_{n+1}}{n}\frac{n(1 - e^{it\log(1+\frac{1}{n})})}{\log(n+1)}\Big|^2+
    \sum^{ }_{n \geq t}  \Big|\frac{c_{n+1}}{n}\frac{n(1 - e^{it\log(1+\frac{1}{n})})}{\log(n+1)}\Big|^2 \\\leq
      E\sum^{ }_{2 \leq n < t}  \Big|\frac{c_{n+1}}{n}\Big|^2\Big(\frac{t}{\log(t)}\Big)^2+  \sum^{ }_{n \geq t}  \Big|\frac{c_{n+1}}{n}\frac{tn\log(1+\frac{1}{n})(1 - e^{it\log(1+\frac{1}{n})})}{\log(n+1)t\log(1+\frac{1}{n})}\Big|^2\\\leq
      E\sum^{ }_{2 \leq n < t}  \Big|\frac{c_{n+1}}{n}\Big|^2\Big(\frac{t}{\log(t)}\Big)^2+  F\sum^{ }_{n \geq t}  \Big|\frac{c_{n+1}}{n}\frac{(1 - e^{it\log(1+\frac{1}{n})})}{t\log(1+\frac{1}{n})}\Big|^2\Big(\frac{t}{\log(t)}\Big)^2\\\leq ( E+G)\sum^{\infty}_{n =  2}  \Big|\frac{c_{n+1}}{n}\Big|^2\Big(\frac{t}{\log(t)}\Big)^2.
\end{gather*}
Where we have used the boundedness of $s\log(1+\frac{1}{s})$ and $\frac{1-e^{is}}{s}$ for $s\in\mathbb{R}^+$. Thus, again due to \eqref{Har2},
\begin{equation*}
\| \widetilde{T}(t) \widetilde{A}^{-1} \widetilde{x}\|_1\leq (B+C+D+4E+4G)^{\frac{1}{2}}\frac{t}{\log(t)}\|\widetilde{x}\|_1.
\end{equation*}
Thus
\begin{equation}\label{right}
\| \widetilde{T}(t) \widetilde{A}^{-1}\|\leq M_0 \frac{t}{\log(t)},
\end{equation}
for some $M_0>0$ and  $t>e$. We will now prove the opposite inequality 
\begin{equation}\label{left}
    m_0\frac{t}{\log(t)}\leq\| \widetilde{T}(t) \widetilde{A}^{-1}\|,
\end{equation}
for some $m_0>0$. First, we observe that due to the reverse triangle inequality, it holds
\begin{gather*}
\| \widetilde{T}(t) \widetilde{A}^{-1} \widetilde{x}\|_1 = \Big( \sum^{\infty}_{n = 2} \Big|c_{n+1}\frac{e^{it\log(n+1)}}{i\log(n+1)}- c_n\frac{e^{it\log(n)}}{i\log(n)}\Big|^2 + |c_2|^2\Big)^{\frac{1}{2}} \\ \geq \Big(\sum^{\infty}_{n = 2}\Big|c_{n+1}\frac{(e^{it\log(n)}\log(n) - e^{it\log(n+1)}\log(n))}{\log(n+1)\log(n)} \Big|^2\Big)^\frac{1}{2}
\\- \Big(\sum^{\infty}_{n = 2} \Big|c_{n+1}\frac{e^{it\log(n+1)}}{i\log(n+1)} -c_{n+1}\frac{e^{it\log(n)}}{i\log(n)}\Big|^2\Big)^\frac{1}{2}  -\Big(\sum^{\infty}_{n = 2} \Big| (c_{n+1}-c_n)\frac{e^{it\log(n)}}{i\log(n)}\Big|^2\Big)^\frac{1}{2}  -  |c_2|.
\end{gather*}
It follows from previous considerations that
\begin{equation*}
    \| \widetilde{T}(t) \widetilde{A}^{-1} \widetilde{x}\|_1\geq \Big(\sum^{\infty}_{n = 2}\Big|c_{n+1}\frac{(e^{it\log(n)}\log(n) - e^{it\log(n+1)}\log(n))}{\log(n+1)\log(n)} \Big|^2\Big)^\frac{1}{2} - C\|\widetilde{x}\|_1,
\end{equation*}
for some $C>0$. Thus, in order to prove \eqref{left}, it suffices to show that 
\begin{equation}\label{random}
    \sum^{\infty}_{n = 2}\Big|c_{n+1}\frac{(e^{it\log(n)} - e^{it\log(n+1)})}{\log(n+1)} \Big|^2\geq  m_1^2\Big(\frac{t}{\log(t)}\Big)^2\|\widetilde{x}\|_1^2,
\end{equation}
for some $m_1>0$ and $t>e$. To this end, we construct for each $t>e$ an element in $\widetilde{H}_1$ in the following way
\begin{gather*}
    \widetilde{x}^{(t)} = (\textgoth{f})\sum^{\infty}_{n = 1} c_n^{(t)} e_n,  \quad j\in\mathbb{N}, \quad \textnormal{ where }\\
\centering
  c^{(t)}_n =
    \begin{cases}
      n & \text{if } n\leq 2t, \\
      4t-n & \text{if }  2t<n \leq 4t,\\
      0 & \text{otherwise.}
    \end{cases}       
\end{gather*}
Observe that
\begin{equation}\label{Xnorm}
\|\widetilde{x}^{(t)}\|^2_1 \leq 4t  .  
\end{equation}
Now, the following estimate holds (cf. \eqref{random})
\begin{gather*}
   \sum^{\infty}_{n = 2}\Big|c^{(t)}_{n+1}\frac{(e^{it\log(n)} - e^{it\log(n+1)})}{\log(n+1)} \Big|^2\geq  \sum_{t\leq n \leq 2t} \Big|t\frac{1 - e^{it\log(1+\frac{1}{n})})}{\log(n+1)}\Big|^2 \\
  \geq   \Big(\frac{t}{\log(4t)}\Big)^2 \sum_{t\leq n \leq 2t}\Big |\frac{1 - e^{it\log(1+\frac{1}{n})})}{it\log(1+\frac{1}{n}) }it\log(1+\frac{1}{n})\Big|^2 \\ \geq \Big(\frac{t}{\log(4t)}\Big)^2 \sum_{t\leq n \leq 2t}\Big |\frac{1 - e^{it\log(1+\frac{1}{n})})}{t\log(1+\frac{1}{n}) }\log(1+\frac{1}{2t})^t\Big|^2 \\ \geq
  \Big(\frac{Ct}{\log(4t)}\Big)^2 \sum_{0\leq n\leq t}\Big |\frac{1 - e^{it\log(1+\frac{1}{n+t})})}{t\log(1+\frac{1}{n+t}) }\Big|^2 \geq \Big(\frac{Ct^2}{\log(4t)}\Big)^2 \sum_{0\leq n \leq t}D \geq \Big(\frac{Ct}{\log(4t)}\Big)^2 \frac{t}{2}D,
\end{gather*}
for $t>e$ and some $C,D>0$ independent of $t>e$. Combining the above with \eqref{random} and \eqref{Xnorm} gives
\begin{equation*}
m_0 \frac{t}{\log(t)}\leq    \frac{\| \widetilde{T}(t) \widetilde{A}^{-1} \widetilde{x}^{(t)}\|_1}{\|\widetilde{x}^{(t)}\|_1}. 
\end{equation*}
For $t>e$. Together with \eqref{right} this shows that
\begin{equation*}\label{equiv}
    m_0 \frac{t}{\log(t)} \leq \| \widetilde{T}(t) \widetilde{A}^{-1}\|\leq M_0 \frac{t}{\log(t)},
\end{equation*}
for $t>e$. This implies, due to \eqref{Hnorm} that
\begin{equation}\label{equiv2}
     m_0' \frac{1}{\log(t)} \leq\frac{ \| \widetilde{T}(t) \widetilde{A}^{-1}\|}{\|\widetilde{T}(t)\|}\
    \leq M_0' \frac{1}{\log(t)},
\end{equation}
for some $m_0',M_0'>0$ and $t>e$ or, equivalently,
\begin{equation*}
    \frac{ \| \widetilde{T}(t) \widetilde{R}_\mu\|}{\|\widetilde{T}(t)\|}\sim \frac{1}{\log(t)},
\end{equation*}
for $t>e$ and arbitrary $\mu \notin \sigma(\widetilde{A})$. Thus 
\begin{equation*}
    \frac{\|\widetilde{T}(t)\widetilde{R}_\mu\|}{\|\widetilde{T}(t)\|}\to 0, \quad \textnormal{ as } \,t\to\infty.
    \end{equation*}

\begin{flushright}
$\square$
\end{flushright}

A similar result can be obtained by the use of Theorem \ref{main}. We are going to check that the semigroup $\widetilde{T}$ meets the assumptions of Theorem~\ref{main} for arbitrary $N \in\mathbb{N}/\{0\} $. Indeed, for each $\lambda_n = i\log n $ one can choose ${\Gamma}_{n}$ surrounding only one point of $\sigma(\widetilde{A})$, namely $\lambda_n$. Note also that,  for $x\in H\subset \widetilde{H}$,
\begin{align*}
    & \widetilde{A}x=Ax,\\
    & \widetilde{R}(\lambda)x=R(\lambda)x,\\
    & \widetilde{P}_{\Gamma_n}x=P_{\Gamma_n}x.
\end{align*}Hence, due to density of $H$ in $\widetilde{H}$
\begin{equation*}
\widetilde{T}(t)\widetilde{P}_{\Gamma_{n}}\widetilde{x}=e^{it\log n}\widetilde{P}_{\Gamma_{n}}\widetilde{x}, \quad \widetilde{x}\in \widetilde{H}.
\end{equation*}
It is easy to see that the function $f(t):=Mt^N+1$ has the properties \eqref{a}, \eqref{b}, has a concave logarithm, and that the following holds:
\begin{equation*}
   \frac{ \| \widetilde{T}(t)\widetilde{P}_{\Gamma_{n}} \|}{f(t)}\leq \frac{\|\widetilde{P}_{\Gamma_n}\|}{f(t)} \leq\frac{\|\widetilde{P}_{\Gamma_n}\|}{Mt^N}\to 0, \quad \textnormal{ as }t\to\infty,\quad n>0.
\end{equation*}
 This means that the semigroup meets the assumption (b) of Theorem  \ref{main}. Application of the presented result yields 
\begin{equation*}
    0= \underset{t \to +\infty}{\lim} \frac{\|\widetilde{T}(t)\widetilde{R}_\mu\|}{Mt^N+1} = \underset{t \to +\infty}{\lim} \frac{\|\widetilde{T}(t)\widetilde{R}_\mu\|}{t^N}=
    \underset{t \to +\infty}{\lim} \frac{\|\widetilde{T}(t)\widetilde{R}_\mu\|}{\|\widetilde{T}(t)\|},
\end{equation*}
for a fixed $\mu \not\in \sigma(\widetilde{A})$.
\begin{flushright}
$\square$
\end{flushright}
\end{exm}
The application of Theorem \ref{main} rendered much shorter calculations for arbitrary $N$ than calculations "by hand" for the simplest case of $N=1$. One can only expect the calculations to become more complicated for larger $N$. Finally, we state the following corollary concerning a sufficient condition for \eqref{Teza} to hold.

\begin{cor}\label{point}
Let $T$ be a $C_0$-semigroup with the generator $A$ and $\omega_0(T)=0. $ If  
\begin{equation*}
\sigma(A) =\sigma_P(A)= \underset{n\in \mathbb{N}}{\bigcup}\lambda_n\textnormal{ },\quad \lambda_n\in (i\mathbb{R} ),\end{equation*}
where all  eigenvalues $\lambda_n$ are simple, and
\begin{equation*}
\|T(t)\|\sim f(t),
\end{equation*}
for an unbounded function $f(t)$ with concave downwards logarithm, (e.g. $f(t) = t^\alpha,\, \alpha\in \mathbb{R}^+$), 
then \eqref{Teza} holds, i.e.,
\begin{equation*}
    \frac{\|T(t)R_\mu\|}{\|T(t)\|}\to 0, \quad \textnormal{ as } t\to\infty.
    \end{equation*}
\end{cor}



%
%
\bibliographystyle{spmpsci}

\end{document}